\documentclass[a4paper,10pt]{article}

\usepackage{authblk}
\usepackage[utf8]{inputenc} 
\usepackage[T1]{fontenc}    
\usepackage{indentfirst}
\usepackage{lmodern}
\usepackage{mathtools}
\usepackage{verbatim}
\usepackage{url}            
\usepackage{booktabs}       
\usepackage{amsfonts}       
\usepackage{amsmath}
\usepackage{amssymb}
\usepackage{amsthm}
\usepackage[a4paper]{geometry} 
\usepackage{enumitem}
\usepackage{xcolor}
\usepackage{mathrsfs}
\usepackage{calligra}
\usepackage{geometry}
\usepackage{fancyhdr}
\usepackage{listings}
\usepackage[linktocpage=true]{hyperref}
\hypersetup{
    colorlinks=true,
    linkcolor=blue,
    citecolor=violet,
    urlcolor=blue
}
\usepackage{graphicx}
\usepackage[makeroom]{cancel}
\usepackage{tikz}
\usepackage[inkscapeformat=png]{svg}
\usepackage[normalem]{ulem}
\usepackage{stmaryrd}

\theoremstyle{plain}
\newtheorem{theorem}{Theorem}[section]

\newtheorem{corollary}[theorem]{Corollary}

\theoremstyle{definition}
\newtheorem{definition}[theorem]{Definition}

\theoremstyle{remark}
\newtheorem{remark}[theorem]{Remark}

\makeatletter
\def\keywords{\xdef\@thefnmark{}\@footnotetext}
\makeatother

\title{The topology of finite and infinite-dimensional Stiefel manifolds and the relation to Grassmannians}
\author{Nizar El Idrissi}

\newcommand{\Addresses}{{
    \bigskip
    \footnotesize
    \textbf{Nizar El Idrissi.}
    \par\nopagebreak Laboratoire : Equations aux dérivées partielles, Algèbre et Géométrie spectrales.
    \par\nopagebreak
    Département de mathématiques, faculté des sciences, université Ibn Tofail, 14000 Kénitra.\par\nopagebreak 
    \textit{E-mail address} : \texttt{nizar.elidrissi@uit.ac.ma}
}}

\begin{document}

\maketitle

\begin{abstract}
    Stiefel manifolds arise naturally as spaces of injective operators and as total spaces of principal bundles over Grassmannians. While their finite-dimensional topology is governed by Bott periodicity, the infinite-dimensional theory exhibits a striking collapse phenomenon stemming from Kuiper's contractibility theorem. In this expository article, we present a unified treatment of finite and infinite-dimensional Stiefel manifolds over real and complex Hilbert spaces, emphasizing three structural principles: \\
    (i) the interpretation of Stiefel manifolds as spaces of injective operators, \\
    (ii) the polar decomposition as a canonical factorization yielding a homotopy decomposition, \\
    (iii) the role or function of Grassmannians as classifying spaces for $\operatorname{GL}_n(\mathbb{F})$ in the stable limit. \\
    In particular, we show that the polar decomposition provides a global homeomorphism
    \[ \operatorname{St}(n,H) \cong \operatorname{St}_{\operatorname{orth}}(n,H) \times \operatorname{P}_n(\mathbb{F}) \]
    valid in arbitrary Hilbert dimension, and that this factorization isolates the entire difference between finite and infinite-dimensional topology. Then, we discuss the implications for the homotopy type of Stiefel manifolds and the relation of Stiefel manifolds to the theory of classifying spaces and characteristic classes.
\end{abstract}



\tableofcontents

\section{Introduction}

The Stiefel manifolds, named after the Swiss mathematician Eduard Stiefel (1909-1978), occupy a central position in algebraic topology, differential geometry, and global analysis because they encode linear data in geometric form.  Originally introduced in Stiefel's 1935 dissertation \cite{Stiefel1936} and independently by Whitney \cite{Whitney1935, Whitney1940}, a Stiefel manifold parametrizes $n$-frames in a vector space and serves as the total space of a principal bundle over a Grassmannian -- which turns out to be universal in the infinite-dimensional case. The topological invariants of Stiefel manifolds -- homotopy, homology, and cohomology groups -- provide the fundamental tools for classifying vector bundles and understanding geometric structures on manifolds \cite{Atiyah1967, Bott1959, Gitler1966,GreenbergHarper1981, Husemoller1994, James1976, Milnor1956, MilnorStasheff1974,Novikov1965, Palais1966, Steenrod1951}. \\

While the classical theory restricts attention to finite-dimensional Euclidean spaces $\mathbb{R}^k$ or $\mathbb{C}^k$, the past five decades have witnessed growing interest in infinite-dimensional Stiefel manifolds. This development was catalyzed by Kuiper's landmark 1965 theorem \cite{Kuiper1965}, which established that the unitary group of an infinite-dimensional Hilbert space is contractible, and more generally that $\operatorname{GL}(H)$ is contractible for infinite-dimensional $H$. This result has profound implications: the classifying space $\operatorname{BGL}(H)$ is trivial, and consequently all infinite-dimensional Hilbert space vector bundles over paracompact bases are trivial \cite{Atiyah1967,Eells1974}. The contrast with the rich nontrivial topology of finite-dimensional Stiefel manifolds  -- nonvanishing homotopy groups -- exemplifies one of the most striking phase transitions in geometric topology. \\

In finite dimension, Stiefel manifolds are highly nontrivial spaces. Their homotopy groups are governed by the deep periodic phenomena discovered by Raoul Bott \cite{Bott1959}. The situation changes dramatically in infinite dimension. Kuiper’s theorem asserts that the unitary group of a separable infinite-dimensional Hilbert space is contractible in the norm topology \cite{Kuiper1965}. As a consequence, infinite-dimensional Grassmannians become classifying spaces, and orthonormal Stiefel manifolds collapse homotopically to a point. At first sight, this suggests that the infinite-dimensional theory should be trivial. However, this impression is misleading: although the orthonormal Stiefel manifold becomes contractible, universal vector bundles of finite rank over the Grassmannian remain nontrivial. The topology has not disappeared -- it has reorganized itself in the form of the intricate finite-rank vector bundles topology over Grassmannians. \\

The purpose of this expository article is to explain, in a unified and structural manner, how these two regimes fit into a single conceptual picture. Rather than treating finite and infinite dimension as separate theories, we show that both arise from the same three organizing principles.

\begin{enumerate}
    \item \textbf{Stiefel manifolds as spaces of injective operators.} \\
          Given a Hilbert space $H$ over $\mathbb{F}=\mathbb{R}$ or $\mathbb{C}$, the Stiefel manifold $\operatorname{St}(n,H)$ is most naturally viewed not merely as the space of $n$-frames, but as the space of injective operators
          \[ \operatorname{St}(n,H)= \{ T \in L(\mathbb{F}^n,H) \; / \; T \text{ is injective} \}. \]
          This operator-theoretic viewpoint is decisive. It identifies $\operatorname{St}(n,H)$ as an open subset of the Banach space $L(\mathbb{F}^n,H)$, thereby placing it in a functional-analytic framework valid in arbitrary Hilbert dimension. Moreover, the canonical right action of $\operatorname{GL}_n(\mathbb{F})$ by precomposition makes the Grassmannian appear naturally as a quotient space.

    \item \textbf{The polar decomposition as a global structural factorization.} \\
          Every injective operator $V:\mathbb{F}^n\to H$ admits a canonical polar factorization $V = U \circ P$, where $U$ is an isometric embedding and $P$ is positive definite. This elementary analytic statement has a strong global topological consequence: it produces a homeomorphism
          \[ \operatorname{St}(n,H) \cong \operatorname{St}_{\operatorname{orth}}(n,H) \times \operatorname{P}_n(\mathbb{F}) \]
          valid for arbitrary Hilbert spaces. \\
          This factorization isolates the entire difference between finite and infinite-dimensional topology. The space $\operatorname{P}_n(\mathbb{F})$ is contractible in all cases; hence the homotopy type of $\operatorname{St}(n,H)$ is completely determined by its orthonormal part. In finite dimension, this orthonormal Stiefel manifold retains the rich structure controlled by Bott periodicity. In infinite dimension, Kuiper’s theorem \cite{Kuiper1965} forces the orthonormal factor to become contractible.

    \item\textbf{Grassmannians as stable classifying spaces.} \\
          The quotient of $\operatorname{St}(n,H)$ by $\operatorname{GL}_n(\mathbb{F})$ yields the Grassmannian $\operatorname{Gr}_n(H)$. In finite dimension, Grassmannians approximate classifying spaces in the stable range $\dim(H) \gg n$, meaning that more and more of their homotopy groups become equal to those of $B\operatorname{GL}_n(\mathbb{F})$. In infinite dimension, Kuiper’s theorem implies that $\operatorname{GL}(H)$ is contractible, so that Grassmannians become genuine classifying spaces for finite-rank vector bundles.
\end{enumerate}

From this perspective, finite and infinite-dimensional theories are not opposed but structurally aligned. What changes between finite and infinite dimension is not the formalism, but the homotopy type of the orthonormal factor: the polar decomposition isolates the precise locus where the phase transition occurs. \\

We develop this perspective in three directions, leading to:
\begin{enumerate}
    \item the connectedness and homotopical properties of Stiefel manifolds;
    \item the principal bundle structure and its homotopy-theoretic classification;
    \item the characteristic classes of the tautological vector bundles over Grassmannians in finite versus infinite dimension \cite{MilnorStasheff1974}.
\end{enumerate}

Our aim is not to reprove classical results, but to reorganize them around these structural principles in order to make transparent the transition from classical Bott-periodic topology to the contractible landscape of Hilbertian geometry.

\section{Stiefel manifolds}

Let $\mathbb{F}$ denote either $\mathbb{R}$ or $\mathbb{C}$, and let $H$ be a Hilbert space over $\mathbb{F}$ with inner product $\langle \cdot, \cdot \rangle$. We consider $H$ to be either finite-dimensional ($\dim H = k < \infty$) or infinite-dimensional separable. All topological considerations refer to the norm topology.

\begin{definition}
    For a positive integer $n \leq \dim H$ (allowing $n < \infty$ even when $\dim H = \infty$), we define:

    \begin{enumerate}
        \item \textbf{The Stiefel manifold of independent $n$-frames:}
              \[
                  \operatorname{St}(n,H) := \{(v_1,\ldots,v_n) \in H^n : \{v_1,\ldots,v_n\} \text{ is linearly independent}\}.
              \]
              This is an open subset of $H^n$, hence a smooth Hilbert manifold when $\dim H = \infty$.

        \item \textbf{The Stiefel manifold of orthonormal $n$-frames:}
              \[
                  \operatorname{St}_{\operatorname{orth}}(n,H) := \{(v_1,\ldots,v_n) \in H^n : \langle v_i, v_j \rangle = \delta_{ij}\}.
              \]
              When $\dim H = k < \infty$, this is the classical compact Stiefel manifold $V_n(\mathbb{F}^k)$. When $\dim H = \infty$, it is an infinite-dimensional Hilbert manifold \cite{Hamilton1982}.

        \item \textbf{The Grassmannian of $n$-dimensional subspaces:}
              \[
                  \operatorname{Gr}_n(H) := \{E \subset H : E \text{ is a closed linear subspace with } \dim E = n\}.
              \]
              This carries a natural topology via the gap metric or equivalently as a quotient of $\operatorname{St}_{\operatorname{orth}}(n,H)$ \cite{Kato1976}.
    \end{enumerate}
\end{definition}

\begin{remark}
    \mbox{}

    \begin{itemize}
        \item In finite dimensions, $\operatorname{St}(n,\mathbb{F}^k)$ is an open subset of $\mathbb{F}^{nk}$, hence a noncompact manifold of dimension $2nk$ (complex) or $nk$ (real). The orthonormal Stiefel manifold $\operatorname{St}_{\operatorname{orth}}(n,\mathbb{F}^k)$ is compact, with dimension $nk - \frac{n(n+1)}{2}$ (real) or $2nk - n^2$ (complex). The classical notation $V_n(\mathbb{F}^k) := \operatorname{St}_{orth}(n,\mathbb{F}^k)$ is standard in the literature \cite{James1976,MilnorStasheff1974}.
        \item When $H$ is infinite-dimensional, $\operatorname{St}_{\operatorname{orth}}(n,H)$ is contractible (see theorem \ref{theorem-H-infinite-dimensional-Storth-contractible} of the present article). This fundamental fact follows from Kuiper's theorem \cite{Kuiper1965} and will be central to our analysis.
    \end{itemize}
\end{remark}

\section{Polar decomposition}

The key technical tool used for studying Stiefel manifolds is the polar decomposition of an $n$-tuple of linearly independent vectors, viewed as an injective linear map $\mathbb{F}^n \to H$.

\begin{theorem}[Polar decomposition homeomorphism]
    For any Hilbert space $H$ (finite or infinite-dimensional, real or complex) and any integer $n$ with $\dim H \geq n$, there is a homeomorphism
    \[
        \Phi: \operatorname{St}(n,H) \xrightarrow{\cong} \operatorname{St}_{\operatorname{orth}}(n,H) \times \operatorname{P}_n(\mathbb{F}),
    \]
    where $\operatorname{P}_n(\mathbb{F})$ denotes the set of positive definite self-adjoint matrices with coefficients in $\mathbb{F}$.
\end{theorem}

\begin{proof}
    Let $(v_1,\ldots,v_n) \in \operatorname{St}(n,H)$. Define $V: \mathbb{F}^n \to H$ by $V(e_i) = v_i$, where $\{e_i\}$ is the standard basis. Then $V$ is injective with closed range (all finite-dimensional subspaces of a Hilbert space are closed). Consider the positive definite operator $P = (V^*V)^{1/2}$ on $\mathbb{F}^n$, where $V^*: H \to \mathbb{F}^n$ is the Hilbert space adjoint. Standard spectral theory guarantees the existence and uniqueness of $P$. \\
    The polar decomposition of $V$ is $V = U \circ P$, where $U: \mathbb{F}^n \to H$ is a partial isometry with initial space $(\ker V)^\perp = \mathbb{F}^n$ and final space $\operatorname{Im} V$. Explicitly, $U = V \circ P^{-1}$. The key observation is that $U$ is an isometry:
    \[
        \langle Ux, Uy \rangle = \langle P^{-1}x, V^*V P^{-1}y \rangle = \langle P^{-1}x, P^2 P^{-1}y \rangle = \langle x, y \rangle.
    \]
    Now interpret:
    \begin{itemize}
        \item $U(e_1),\ldots,U(e_n)$ form an orthonormal $n$-frame in $H$, hence belong to $\operatorname{St}_{\operatorname{orth}}(n,H)$;
        \item $P$, expressed in the standard basis, is a positive definite $n \times n$ matrix.
    \end{itemize}
    The map $\Phi(V) = (U,P)$ is clearly continuous in both directions (matrix inversion and square root are continuous on positive definite matrices), with continuous inverse $(U,P) \mapsto U \circ P$.
\end{proof}

\begin{remark}
    This theorem is well-known in finite dimensions, but the proof given above works without modification when $H$ is infinite-dimensional. The point is that $V^*V$ remains a positive definite and invertible operator on $\mathbb{F}^n$; moreover, its square root is defined via the functional calculus and depends continuously on $V$.
\end{remark}

\begin{corollary}
    For any Hilbert space $H$ (finite or infinite-dimensional, real or complex) and any integer $n$ with $\dim H \geq n$, the homotopy type of $\operatorname{St}(n,H)$ is the product of the homotopy types of $\operatorname{St}_{\operatorname{orth}}(n,H)$ and $\operatorname{P}_n(\mathbb{F})$. Consequently, $\operatorname{St}(n,H)$ is homotopy equivalent to $\operatorname{St}_{\operatorname{orth}}(n,H)$ since $\operatorname{P}_n(\mathbb{F})$ is convex.
\end{corollary}

\section{Connectedness and homotopical properties of orthonormal Stiefel manifolds}

The connectedness properties of $\operatorname{St}_{\operatorname{orth}}(n,H)$ are fundamental and exhibit a sharp transition depending on whether $\dim H \geq n+1$ (real case) or $\dim H \geq n$ (complex case).

\begin{theorem}[Connectedness of real orthonormal Stiefel manifolds]
    For $\mathbb{F} = \mathbb{R}$, $\operatorname{St}_{\operatorname{orth}}(n,\mathbb{R}^k)$ is connected if and only if $k \geq n+1$. For $k = n$, it consists of two connected components distinguished by the sign of the determinant of the $n \times n$ matrix of frame vectors.
\end{theorem}

\begin{proof}
    The group $O(k)$ acts transitively on $\operatorname{St}_{\operatorname{orth}}(n,\mathbb{R}^k)$ with stabilizer $O(k-n)$ (embedded as matrices fixing the first $n$ standard basis vectors). Hence,
    \[
        \operatorname{St}_{\operatorname{orth}}(n,\mathbb{R}^k) \cong O(k)/O(k-n).
    \]
    Consider the exact homotopy sequence of the fibration $O(k-n) \to O(k) \to \operatorname{St}_{\operatorname{orth}}(n,\mathbb{R}^k)$:
    \[
        \pi_0(O(k-n)) \to \pi_0(O(k)) \to \pi_0(\operatorname{St}_{\operatorname{orth}}(n,\mathbb{R}^k)) \to 0.
    \]
    We know $\pi_0(O(m)) \cong \mathbb{Z}_2$ for all $m \geq 1$, with the nontrivial element represented by a reflection. The map $\pi_0(O(k-n)) \to \pi_0(O(k))$ is induced by inclusion. This map is:
    \begin{itemize}
        \item An isomorphism when $k-n \geq 1$ (since both groups are $\mathbb{Z}_2$ and the inclusion sends the nontrivial component to the nontrivial component);
        \item The trivial map when $k-n = 0$ (i.e., $k=n$), because $O(0)$ is a point with trivial $\pi_0$.
    \end{itemize}
    Exactness then gives:
    \begin{itemize}
        \item If $k-n \geq 1$, the sequence is $\mathbb{Z}_2 \xrightarrow{\cong} \mathbb{Z}_2 \to \pi_0(\operatorname{St}_{\operatorname{orth}}) \to 0$, forcing $\pi_0(\operatorname{St}_{\operatorname{orth}}) = 0$;
        \item If $k=n$, we have $0 \to \mathbb{Z}_2 \to \pi_0(\operatorname{St}_{\operatorname{orth}}) \to 0$, so $\pi_0(\operatorname{St}_{\operatorname{orth}}) \cong \mathbb{Z}_2$.
    \end{itemize}
    Thus, $\operatorname{St}_{\operatorname{orth}}(n,\mathbb{R}^n) \cong O(n)$ consists of two components (matrices with determinant $\pm 1$), while for $k \geq n+1$ it is connected.
\end{proof}

\begin{theorem}[Connectedness of complex orthonormal Stiefel manifolds]
    For $\mathbb{F} = \mathbb{C}$, $\operatorname{St}_{\operatorname{orth}}(n,\mathbb{C}^k)$ is connected for all $k \geq n$. In particular, $U(n)$ is connected.
\end{theorem}

\begin{proof}
    The unitary group $U(k)$ is connected for all $k$. The transitive action $U(k) \to \operatorname{St}_{\operatorname{orth}}(n,\mathbb{C}^k)$ with stabilizer $U(k-n)$ gives a fibration $U(k-n) \to U(k) \to \operatorname{St}_{\operatorname{orth}}(n,\mathbb{C}^k)$. Since $\pi_0(U(k)) = 0$ and $\pi_0(U(k-n)) = 0$, the exact sequence yields $\pi_0(\operatorname{St}_{\operatorname{orth}}) = 0$.
\end{proof}

In fact, the homotopy groups of Stiefel manifolds are classical and completely known.

\begin{theorem}[Homotopy groups of Stiefel manifolds]
    For $k \geq n$, the following hold:

    \begin{enumerate}
        \item \textbf{Real case:} $\pi_i(\operatorname{St}_{orth}(n,\mathbb{R}^k)) = 0$ for $i \leq k-n-1$, and $\pi_{k-n}(\operatorname{St}_{orth}(n,\mathbb{R}^k)) \cong \mathbb{Z}$ if $k-n$ is even or $n=1$; otherwise $\pi_{k-n}(\operatorname{St}_{orth}(n,\mathbb{R}^k)) \cong \mathbb{Z}_2$. More precisely,
              \[
                  \pi_{k-n}(\operatorname{St}_{orth}(n,\mathbb{R}^k)) \cong
                  \begin{cases}
                      \mathbb{Z}   & \text{if } n=1 \text{ or } k-n \text{ is even},      \\[2mm]
                      \mathbb{Z}_2 & \text{if } n \geq 2 \text{ and } k-n \text{ is odd}.
                  \end{cases}
              \]
        \item \textbf{Complex case:} $\pi_i(\operatorname{St}_{orth}(n,\mathbb{C}^k)) = 0$ for $i \leq 2(k-n)$, and $\pi_{2(k-n)+1}(\operatorname{St}_{orth}(n,\mathbb{C}^k)) \cong \mathbb{Z}$.
    \end{enumerate}
\end{theorem}

\begin{proof}
    These results follow from the homotopy exact sequence of the fibration $O(k-n) \to O(k) \to \operatorname{St}_{orth}(n,\mathbb{R}^k)$ and the stable homotopy groups of orthogonal groups (Bott periodicity). For the complex case, use $U(k-n) \to U(k) \to \operatorname{St}_{orth}(n,\mathbb{C}^k)$ and Bott periodicity: $\pi_i(U) \cong \mathbb{Z}$ for odd $i$, $0$ for even $i$ \cite{Bott1959, Palais1966}.
\end{proof}

In the infinite-dimensional case, we have

\begin{theorem}[Infinite-dimensional case]
    \label{theorem-H-infinite-dimensional-Storth-contractible}
    For $\dim H = \infty$, $\operatorname{St}_{\operatorname{orth}}(n,H)$ is contractible, hence connected, for both real and complex scalars.
\end{theorem}

\begin{proof}
    The group of orthogonal (respectively unitary) operators on a separable infinite-dimensional Hilbert space is contractible in the norm topology \cite{Kuiper1965}. This group acts transitively on $\operatorname{St}_{\operatorname{orth}}(n,H)$ with stabilizer isomorphic to the same group (since $H^\perp \cong H$). Hence, $\operatorname{St}_{\operatorname{orth}}(n,H)$ is a homogeneous space of a contractible group, implying it is also contractible. More directly, one can construct an explicit contraction using Kuiper's theorem or the Eilenberg swindle \cite{Eells1974,Kuiper1965}.
\end{proof}

\section{Principal bundle structure}

Let $n$ be a natural number and $H$ be a Hilbert space with $\dim(H) \geq n$. The topology on $\mathrm{Gr}_n(H)$ is induced by the gap metric
\[
  \delta(E, F) \;=\; \|P_E - P_F\|,
\]
where $P_E$ (resp.\ $P_F$) denotes the orthogonal projection onto $E$ (resp.\ $F$),
and $\|\cdot\|$ is the operator norm on $\mathcal{L}(H)$. Equivalently, this
topology coincides with the quotient topology arising from the surjection
$\mathrm{St}_{\mathrm{orth}}(n, H) \to \mathrm{Gr}_n(H)$; we use both
descriptions freely below.

We now establish the principal bundle structure.

\begin{theorem}[Principal bundle structure]
    For any Hilbert space $H$ (any dimension, $\mathbb{F} = \mathbb{R}$ or $\mathbb{C}$) and $\dim H \geq n$, the map
\[
  \pi \colon \mathrm{St}(n, H) \longrightarrow \mathrm{Gr}_n(H), \qquad
  \pi(v_1, \ldots, v_n) \;=\; \mathrm{span}\{v_1, \ldots, v_n\},
\]
is a locally trivial principal $\mathrm{GL}(n, \mathbb{F})$-bundle, where
$\mathrm{GL}(n, \mathbb{F})$ acts on the right by
\[
  (v_1, \ldots, v_n) \cdot g \;=\;
  \Bigl(\sum_{j=1}^n g_{j1}\,v_j,\; \ldots,\; \sum_{j=1}^n g_{jn}\,v_j\Bigr),
  \qquad g = (g_{ij}) \in \mathrm{GL}(n, \mathbb{F}).
\]
In operator notation: $V \cdot g = V \circ g$ for $V \in \mathrm{St}(n, H)$ and
$g \in \mathrm{GL}(n, \mathbb{F})$.
\end{theorem}

\begin{proof}
    We verify the four defining properties of a locally trivial principal bundle
in turn: continuity of $\pi$, freeness and transitivity of the
$\mathrm{GL}(n, \mathbb{F})$-action on fibers, and local triviality.

\begin{itemize} 
    \item \textbf{Step 1. Continuity of $\pi$}

We show that $\pi$ is continuous with respect to the norm topology on
$\mathrm{St}(n, H) \subset H^n$ and the gap metric on $\mathrm{Gr}_n(H)$.

Let $V_\alpha \to V$ in $\mathrm{St}(n, H)$, i.e., $V_\alpha(e_i) \to V(e_i)$ in $H$
for each $i \in \{1, \ldots, n\}$. Write $E_\alpha = \mathrm{Im}\, V_\alpha$ and
$E = \mathrm{Im}\, V$. We claim that $\delta(E_\alpha, E) \to 0$, i.e.,
$P_{E_\alpha} \to P_E$ in operator norm.

For any injective operator $T \colon \mathbb{F}^n \to H$ with finite-dimensional
image, the orthogonal projection onto $\mathrm{Im}\, T$ is given by
\[
  P_{\mathrm{Im}\, T} = T(T^*T)^{-1}T^*.
\]
The map $T \mapsto T(T^*T)^{-1}T^*$ is continuous on the open set
$\mathcal{L}(\mathbb{F}^n, H)_{\mathrm{inj}}$ of injective operators: indeed,
$T^* T$ is a positive definite $n \times n$ matrix depending continuously on $T$
(the adjoint map and operator composition are both norm-continuous), matrix
inversion is continuous on $\mathrm{GL}(n, \mathbb{F})$, and composition of
bounded operators is jointly norm-continuous. Therefore,
\[
  P_{E_\alpha}
  = V_\alpha (V_\alpha^* V_\alpha)^{-1} V_\alpha^*
  \;\longrightarrow\;
  V(V^*V)^{-1}V^*
  = P_E
\]
in operator norm as $V_\alpha \to V$, which gives $\delta(E_\alpha, E) \to 0$.
Hence $\pi$ is continuous.

\item \textbf{Step 2. Right action, freeness, and transitivity on fibers}

\medskip
\noindent\textit{Well-definedness of the action.}
For $V \in \mathrm{St}(n, H)$ and $g \in \mathrm{GL}(n, \mathbb{F})$, set
$V \cdot g := V \circ g$. Since $g$ is an automorphism of $\mathbb{F}^n$ and $V$
is injective, the composition $V \circ g$ is again injective, so
$V \cdot g \in \mathrm{St}(n, H)$. The map $(V, g) \mapsto V \circ g$ is
bilinear and jointly continuous in the norm topology. The axioms
$(V \cdot g) \cdot h = V \cdot (gh)$ and $V \cdot \mathrm{id} = V$ hold by
associativity of composition.

\medskip
\noindent\textit{The action preserves fibers.}
For any $g \in \mathrm{GL}(n, \mathbb{F})$,
\[
  \mathrm{Im}(V \circ g) = V\!\left(g(\mathbb{F}^n)\right) = V(\mathbb{F}^n)
  = \mathrm{Im}\, V,
\]
since $g$ is surjective. Hence $\pi(V \cdot g) = \pi(V)$, and the action
restricts to each fiber $\pi^{-1}(E)$.

\medskip
\noindent\textit{Freeness.}
Suppose $V \cdot g = V$, i.e., $V \circ g = V$. Then for each basis vector $e_i$,
\[
  V(g(e_i)) = V(e_i).
\]
Injectivity of $V$ gives $g(e_i) = e_i$ for all $i = 1, \ldots, n$, hence
$g = \mathrm{id}_{\mathbb{F}^n}$. The action is therefore free.

\medskip
\noindent\textit{Transitivity on fibers.}
Fix $E \in \mathrm{Gr}_n(H)$ and let $V, W \in \pi^{-1}(E)$, so that
$\mathrm{Im}\, V = \mathrm{Im}\, W = E$. Both $V$ and $W$ are linear isomorphisms
$\mathbb{F}^n \xrightarrow{\;\sim\;} E$. Set
\[
  g \;:=\; V^{-1} \circ W \;\in\; \mathrm{GL}(n, \mathbb{F}),
\]
where $V^{-1} \colon E \to \mathbb{F}^n$ is the inverse of $V$ viewed as an
isomorphism onto $E$. Then
\[
  V \cdot g = V \circ (V^{-1} \circ W) = W.
\]
Uniqueness of $g$ follows from freeness. Hence the action is simply transitive on
each fiber, identifying $\pi^{-1}(E) \cong \mathrm{GL}(n, \mathbb{F})$ as a
$\mathrm{GL}(n, \mathbb{F})$-torsor.

 \item \textbf{Step 3. Local triviality}

This is the heart of the proof. We construct an explicit local trivialization
around an arbitrary point $E_0 \in \mathrm{Gr}_n(H)$, compatible with the gap
metric topology.

\medskip
\noindent\textbf{Construction of the trivializing neighborhood.}
Let $E_0 \in \mathrm{Gr}_n(H)$, decompose $H = E_0 \oplus E_0^\perp$ orthogonally. Define
$$
\Omega := \bigl\{\, E \in \mathrm{Gr}_n(H) \;\big|\; P\big|_E : E \to E_0 \text{ is an isomorphism} \,\bigr\}.
$$

Since $(P_{E_0})|_{E_0} = \mathrm{id}_{E_0}$, we have $E_0 \in \Omega$. We show that $\Omega$ is open in the gap metric.

Let $E \in \Omega$. Since $\dim E = n < \infty$ and $P|_E : E \to E_0$ is an isomorphism between $n$-dimensional spaces, it is in particular bounded below: there exists $c > 0$ such that
$$
\|P x\| \geq c\,\|x\| \qquad \text{for all } x \in E \qquad (*)
$$

Now let $F \in \mathrm{Gr}_n(H)$ with $\delta(E, F) = \|P_E - P_F\| < \varepsilon$, for $\varepsilon > 0$ to be chosen. Take any $y \in F$ with $\|y\| = 1$. We estimate $\|Py\|$ from below in two steps.

\textbf{First step.} Since $y \in F$, we have $P_F y = y$, so
$$
\|P_E y\| = \|P_E y - P_F y + y\| \geq \|y\| - \|(P_E - P_F)y\| \geq 1 - \delta(E, F).
$$
In particular, $P_E y \in E$ with $\|P_E y\| \geq 1 - \delta(E, F)$.

\textbf{Second step.} Apply the bounded-below estimate $(*)$ to $P_E y \in E$:
$$
\|P(P_E y)\| \geq c\,\|P_E y\| \geq c\,(1 - \delta(E, F)).
$$
Then, since $\|P\| \leq 1$,
$$
\|Py\| \geq \|P(P_E y)\| - \|P(y - P_E y)\|
\geq c\,(1 - \delta(E,F)) - \|y - P_E y\|.
$$
Again using $P_F y = y$:
$$
\|y - P_E y\| = \|(P_F - P_E)y\| \leq \delta(E, F).
$$
Combining, we obtain
$$
\|Py\| \geq c - (c + 1)\,\delta(E, F).
$$

Choosing $\varepsilon := \dfrac{c}{2(c+1)}$, for any $F$ with $\delta(E, F) < \varepsilon$ and any unit $y \in F$:
$$
\|Py\| \geq \frac{c}{2} > 0.
$$
Hence $P|_F$ is bounded below by $c/2$ on the $n$-dimensional space $F$, so it is injective. Since $\dim F = \dim E_0 = n$, it is an isomorphism $F \to E_0$, giving $F \in \Omega$. Therefore, $\Omega$ is open.

\medskip
\noindent\textbf{Construction of a continuous section.}
Fix an ordered basis $(b_1, \ldots, b_n)$ of $E_0$. For each $E \in \Omega$,
the map $(P_{E_0})|_E \colon E \to E_0$ is an isomorphism, so define
\[
  s(E) \;:=\; \bigl(((P_{E_0})|_E)^{-1}b_1,\; \ldots,\; ((P_{E_0})|_E)^{-1}b_n\bigr)
  \;\in\; H^n,
\]
where $((P_{E_0})|_E)^{-1} \colon E_0 \to E$ is the inverse isomorphism.

We verify that $s(E) \in \pi^{-1}(E) \subset \mathrm{St}(n, H)$:
\begin{enumerate}[label=(\roman*)]
\item Each vector $((P_{E_0})|_E)^{-1}b_i$ lies in $E$ by definition of $((P_{E_0})|_E)^{-1}$.
\item The vectors $((P_{E_0})|_E)^{-1}b_1, \ldots, ((P_{E_0})|_E)^{-1}b_n$ are linearly
  independent: if $\sum_i \lambda_i ((P_{E_0})|_E)^{-1}b_i = 0$, applying $(P_{E_0})|_E$
  gives $\sum_i \lambda_i b_i = 0$, hence all $\lambda_i = 0$ since
  $(b_i)$ is a basis of $E_0$.
\item Their span equals $E$: we have $n$ linearly independent vectors in the
  $n$-dimensional space $E$, so they span $E$.
\end{enumerate}
Hence $\pi(s(E)) = E$, confirming that $s$ is a section of $\pi$ over $\Omega$.

Continuity of $s$: the map $E \mapsto (P_{E_0})|_E$ is continuous from
$(\Omega, \delta)$ to $\mathcal{L}(E_0, H)$ (since $E \mapsto P_E$ is
norm-continuous and $(P_{E_0})|_E = ((P_E)|_{E_0})^{-1}$), and hence $s$ is continuous.

\medskip
\noindent\textbf{Construction of the local trivialization.}
Define
\[
  \Phi \colon \pi^{-1}(\Omega) \;\longrightarrow\; \Omega \times \mathrm{GL}(n, \mathbb{F})
\]
as follows. Given $V \in \pi^{-1}(\Omega)$, let $E := \pi(V) = \mathrm{Im}\, V
\in \Omega$. Since the action is simply transitive on $\pi^{-1}(E)$ and
$s(E) \in \pi^{-1}(E)$, there exists a unique $g_V \in \mathrm{GL}(n, \mathbb{F})$
such that $V = s(E) \cdot g_V = s(E) \circ g_V$. Set
\[
  \Phi(V) \;:=\; \bigl(\pi(V),\; g_V\bigr).
\]
Explicitly, $g_V = s(\pi(V))^{-1} \circ V$, where $s(\pi(V))^{-1}$ denotes
the inverse of $s(\pi(V)) \colon \mathbb{F}^n \to E$ as an isomorphism. In
coordinates,
\[
  (g_V)_{ij} \;=\; e_i^*\!\left(s(\pi(V))^{-1}(V(e_j))\right),
\]
a composition of continuous maps in $V$, so $\Phi$ is continuous.

The inverse map is
\[
  \Psi \colon \Omega \times \mathrm{GL}(n, \mathbb{F}) \;\longrightarrow\;
  \pi^{-1}(\Omega), \qquad
  \Psi(E, g) \;=\; s(E) \circ g,
\]
which is continuous since $s$ is continuous and operator composition is jointly
continuous. One verifies directly that $\Phi \circ \Psi = \mathrm{id}$ and
$\Psi \circ \Phi = \mathrm{id}$, so $\Phi$ is a homeomorphism.

\medskip
\noindent\textbf{Equivariance of $\Phi$.}
For $h \in \mathrm{GL}(n, \mathbb{F})$ and $V \in \pi^{-1}(\Omega)$:
\[
  \Phi(V \cdot h)
  = \bigl(\pi(V \cdot h),\; g_{V \cdot h}\bigr)
  = \bigl(\pi(V),\; g_V \cdot h\bigr),
\]
since $s(\pi(V)) \circ (g_V h) = (s(\pi(V)) \circ g_V) \circ h = V \circ h =
V \cdot h$, confirming uniqueness of $g_{V \cdot h} = g_V h$. Thus $\Phi$
intertwines the right $\mathrm{GL}(n, \mathbb{F})$-action on $\pi^{-1}(\Omega)$
with the natural right action on the second factor of
$\Omega \times \mathrm{GL}(n, \mathbb{F})$.

\medskip
\noindent\textbf{Covering of $\mathrm{Gr}_n(H)$.}
For each $E_0 \in \mathrm{Gr}_n(H)$, the above construction yields an open
neighborhood $\Omega_{E_0}$ of $E_0$ over which the bundle trivializes. The
collection $\{\Omega_{E_0}\}_{E_0 \in \mathrm{Gr}_n(H)}$ covers
$\mathrm{Gr}_n(H)$.
\end{itemize}

\end{proof}

\begin{corollary}[Homotopy classification of the bundle]
    The principal $\operatorname{GL}(n,\mathbb{F})$-bundle $\operatorname{St}(n,\mathbb{F}^k) \to \operatorname{Gr}_n(\mathbb{F}^k)$ is classified by a map $f: \operatorname{Gr}_n(\mathbb{F}^k) \to B\operatorname{GL}(n,\mathbb{F})$. The induced map on homotopy groups in the stable range $k \gg n$ is an isomorphism. In particular:

    \begin{itemize}
        \item For $\mathbb{F} = \mathbb{R}$, the classifying map is an isomorphism on $\pi_i$ for $i \leq k-n-1$;
        \item For $\mathbb{F} = \mathbb{C}$, the classifying map is an isomorphism on $\pi_i$ for $i \leq 2(k-n)$.
    \end{itemize}
    This follows from the fact that $\operatorname{St}_{\operatorname{orth}}(n,\mathbb{F}^k)$ is $(k-n-1)$-connected (real) or $(2k-2n)$-connected (complex) \cite{James1976}.
\end{corollary}

\begin{remark}
  In the stable range $k \gg n$, the Grassmannian approximates the classifying space, meaning that more and more of its homotopy groups become equal to those of $B\operatorname{GL}_n(\mathbb{F})$.
\end{remark}

\section{Characteristic classes of the tautological vector bundles over Grassmannians}

For any Hilbert space $H$ over $\mathbb{R}$ or $\mathbb{C}$ and any positive integer $n$ with $\dim H \geq n$, let $\operatorname{Gr}_n(H)$ be the Grassmannian of $n$-dimensional subspaces of $H$. The tautological vector bundle $\gamma^n \to \operatorname{Gr}_n(H)$ is defined such that the fiber over each subspace $V \in \operatorname{Gr}_n(H)$ is $V$ itself. \\

We distinguish between the cases where $H$ is infinite-dimensional and where $H$ is finite-dimensional, say $\dim H = k \geq n$. In the infinite-dimensional case, $\operatorname{Gr}_n(H)$ serves as the classifying space for vector bundles, leading to a polynomial cohomology ring generated by the characteristic classes without relations. In the finite-dimensional case, denoted $\operatorname{Gr}(n,k)$, the cohomology ring is a quotient of this polynomial ring, imposing relations due to the finite dimensionality.

\begin{enumerate}

    \item \textbf{Infinite-dimensional case}. \\
          When $H$ is infinite-dimensional, $\operatorname{Gr}_n(H)$ is homotopy equivalent to the classifying space $BU(n)$ in the complex case and $BO(n)$ in the real case. The tautological bundle $\gamma^n$ is the universal vector bundle.

          \begin{itemize}
              \item \textbf{Complex case}. \\
                    The cohomology ring is
                    \[
                        H^*(\operatorname{Gr}_n(H); \mathbb{Z}) \cong \mathbb{Z}[c_1, \dots, c_n],
                    \]
                    where $c_i = c_i(\gamma^n) \in H^{2i}(Gr_n(H); \mathbb{Z})$ are the Chern classes, which are algebraically independent. The Chern classes vanish for $i > n$, and the total Chern class is
                    \[
                        c(\gamma^n) = 1 + c_1 + \dots + c_n.
                    \]

              \item \textbf{Real case}. \\
                    The mod-2 cohomology ring is
                    \[
                        H^*(\operatorname{Gr}_n(H); \mathbb{Z}/2\mathbb{Z}) \cong \mathbb{Z}/2\mathbb{Z}[w_1, \dots, w_n],
                    \]
                    where $w_i = w_i(\gamma^n) \in H^i(Gr_n(H); \mathbb{Z}/2\mathbb{Z})$ are the Stiefel-Whitney classes, which are algebraically independent. These vanish for $i > n$, and the total Stiefel-Whitney class is
                    \[
                        w(\gamma^n) = 1 + w_1 + \dots + w_n.
                    \]
                    The rational cohomology is generated by the Pontryagin classes $p_i(\gamma^n) \in H^{4i}(\operatorname{Gr}_n(H); \mathbb{Z})$ for $i = 1, \dots, \lfloor n/2 \rfloor$, with the total Pontryagin class
                    \[
                        p(\gamma^n) = 1 + p_1 + \dots + p_{\lfloor n/2 \rfloor},
                    \]
                    and $p_i(\gamma^n) = 0$ for $i > \lfloor n/2 \rfloor$. When $n$ is even, the Euler class $e(\gamma^n) \in H^n(Gr_n(H); \mathbb{Z})$ also contributes, subject to 2-torsion relations in the integral cohomology.
          \end{itemize}

    \item \textbf{Finite-dimensional case} \\
          When $\dim H = k < \infty$, $\operatorname{Gr}_n(H) = \operatorname{Gr}(n,k)$ is compact, and the tautological bundle $\gamma^n$ satisfies $\gamma^n \oplus Q \cong \underline{\mathbb{F}^k}$, where $Q$ is the quotient bundle of rank $k-n$ (with $\mathbb{F} = \mathbb{R}$ or $\mathbb{C}$). This imposes relations on the characteristic classes, making the cohomology ring a quotient of the infinite-dimensional case. Unlike the infinite case, $\operatorname{Gr}(n,k)$ is not universal for all vector bundles but classifies those over spaces of sufficiently low dimension.

          \begin{itemize}
              \item \textbf{Complex case}. \\
                    The cohomology ring is
                    \[
                        H^*(\operatorname{Gr}(n,k); \mathbb{Z}) \cong \mathbb{Z}[c_1, \dots, c_n] / I,
                    \]
                    where $c_i = c_i(\gamma^n)$, and $I$ is the ideal generated by the relations from $c(\gamma^n) \cdot c(Q) = 1$. The relations are explicitly given by
                    \[
                        \sum_{j=0}^i (-1)^j c_j(\gamma^n) \, c_{k-n+i-j}(Q) = 0
                    \]
                    for appropriate degrees. The ring is graded with dimension $\dim_\mathbb{C} \operatorname{Gr}(n,k) = n(k-n)$, so classes vanish above degree $2n(k-n)$. The generators $c_i$ are not algebraically independent in general, and the relations reflect the fact that $\operatorname{Gr}(n,k)$ only classifies bundles over spaces of dimension at most $2n(k-n)$.

              \item \textbf{Real case}. \\
                    The mod-2 cohomology is
                    \[
                        H^*(\operatorname{Gr}(n,k); \mathbb{Z}/2\mathbb{Z}) \cong \mathbb{Z}/2\mathbb{Z}[w_1, \dots, w_n] / J,
                    \]
                    where $w_i = w_i(\gamma^n)$, and $J$ includes relations from $w(\gamma^n) \cdot w(Q) = 1$, with $w_i = 0$ for $i > n$. The dimension is $\dim_\mathbb{R} Gr(n,k) = n(k-n)$, so classes vanish above this degree. For integral cohomology, Pontryagin classes $p_i(\gamma^n)$ generate the torsion-free part, with relations analogous to the complex case, including $p(\gamma^n) \cdot p(Q) = 1$ and potential 2-torsion. The Euler class appears when orientable, but with finite-dimensional constraints.
          \end{itemize}

    \item \textbf{Differences between infinite and finite cases}. \\
          In the infinite-dimensional case, the cohomology rings are free polynomial rings on the characteristic classes, reflecting the universal classifying property: every vector bundle is pulled back from $\gamma^n$. No relations bind the classes beyond their natural degrees. In contrast, the finite-dimensional case introduces relations via the triviality of the ambient bundle, leading to a quotient ring structure. For large $k$, $\operatorname{Gr}(n,k)$ approximates the infinite case, with relations becoming negligible in low degrees.
\end{enumerate}

\section{Nontriviality of the principal bundle $\operatorname{St}(n,H) \to \operatorname{Gr}_n(H)$}

\begin{theorem}
    For any Hilbert space with $\dim(H) \geq n \geq 1$, the principal $\operatorname{GL}(n,\mathbb{F})$-bundle $\operatorname{St}(n,H) \to \operatorname{Gr}_n(H)$ is trivializable if and only if $n = \dim(H)$.
\end{theorem}

\begin{proof}
    When $  n = k  $, $  \operatorname{Gr}_n(\mathbb{F}^n)  $ is a single point (the entire space $  \mathbb{F}^n  $). The total space $  \operatorname{St}(n,\mathbb{F}^n)  $ is the space of bases for $  \mathbb{F}^n  $, homeomorphic to $  \mathrm{GL}(n,\mathbb{F})  $. The bundle reduces to $  \mathrm{GL}(n,\mathbb{F}) \to \{\mathrm{pt}\}  $, which is trivially the product bundle over a point. \\
    Now, suppose that $H$ is a real Hilbert space with $\dim(H) > n$ and the principal bundle is trivializable. The total space $  St(n,H)  $ is connected (as it retracts after deformation to the compact Stiefel manifold, which is connected when $\dim(H) > n$). The base $  Gr_n(H)  $ is connected, and the fiber $  \mathrm{GL}(n,\mathbb{R})  $ has two components. If the bundle were trivializable, the total space would thus have two components, a contradiction. Hence, the bundle is not trivializable for a real Hilbert space with $\dim(H) > n  $. \\
    Finally, suppose that $H$ is a complex Hilbert space with $\dim(H) > n$. Since the tautological bundle has nonzero Chern classes $c_i(\gamma^n)$ generating the cohomology ring, the bundle cannot be trivializable.
\end{proof}

\section{Examples}

This section illustrates the theory with concrete examples spanning the finite-dimensional, infinite-dimensional, real, and complex cases. Unless $\dim(H)=n$, the bundle $\operatorname{St}(n,H) \to \operatorname{Gr}_n(H)$ always remains nontrivial.

\begin{enumerate}
    \item \textbf{$n=k$, real or complex cases}. \\
          When $\dim H = n$, $\operatorname{St}(n,\mathbb{F}^n) = \operatorname{GL}(n,\mathbb{F})$, $\operatorname{St}_{\operatorname{orth}}(n,\mathbb{F}^n) = O(n,\mathbb{F})$, and $\operatorname{Gr}_n(\mathbb{F}^n)$ is a single point. The bundle collapses to $\operatorname{GL}(n,\mathbb{F}) \to \operatorname{pt}$, which is a trivial bundle.
          The classical polar decomposition writes as:
          \[ \operatorname{GL}(n,\mathbb{F}) \cong O(n,\mathbb{F}) \times \operatorname{P}_n(\mathbb{R}). \]

    \item \textbf{$n=1$, $k>1$, real case}. \\
          When $n=1$ and $k>1$, $\operatorname{St}(1,\mathbb{R}^k) = \mathbb{R}^k \setminus \{0\} \simeq S^{k-1}$, $\operatorname{St}_{\operatorname{orth}}(1,\mathbb{R}^k) = S^{k-1}$, and $\operatorname{Gr}_1(\mathbb{R}^k) = \mathbb{RP}^{k-1}$. The principal bundle $\mathbb{R}^k \setminus \{0\} \to \mathbb{RP}^{k-1}$ with fiber group $\operatorname{GL}(1,\mathbb{R}) \simeq S^0$ is nontrivial since $\mathbb{R}^k \setminus \{0\}$ is connected ($k >1$) and $S^0$ has two components.
          The polar decomposition writes as:
          \[
              \mathbb{R}^k \setminus \{0\} \cong S^{k-1} \times \mathbb{R}_{>0}.
          \]

    \item \textbf{$n=1$, $k>1$, complex case}. \\
          For $n=1$ and $k>1$, $\operatorname{St}(1,\mathbb{C}^k) = \mathbb{C}^k \setminus \{0\} \simeq S^{2k-1}$, $\operatorname{St}_{\operatorname{orth}}(1,\mathbb{C}^k) = S^{2k-1}$, and $\operatorname{Gr}_1(\mathbb{C}^k) = \mathbb{CP}^{k-1}$. The bundle $\mathbb{C}^k \setminus \{0\} \to \mathbb{CP}^{k-1}$ with fiber group $\operatorname{GL}(1,\mathbb{C}) = \mathbb{C}^\times \simeq S^1$ is nontrivial since its Chern class $c_1$ is the generator of $H^2(\mathbb{CP}^{k-1};\mathbb{Z})$. The polar decomposition writes as:
          \[
              \mathbb{C}^k \setminus \{0\} \cong S^{2k-1} \times \mathbb{R}_{>0}.
          \]

    \item \textbf{Infinite-dimensional real or complex case}. \\
          Let $H$ be a separable Hilbert space. The polar decomposition implies that $\operatorname{St}(n,H)$ is contractible. Thus, the bundle $\operatorname{St}(n,H) \to \operatorname{Gr}_n(H)$ is the universal principal bundle with structure group $\operatorname{GL}(n,\mathbb{F})$. It is not trivializable.

          \begin{remark}
              Kuiper's theorem does not imply triviality of this bundle. What Kuiper's theorem implies is that $\operatorname{GL}(H)$ is contractible, which implies that the stabilized bundle of the tautological vector bundle over $\operatorname{Gr}_n(H)$ (adding a complementary infinite-dimensional trivial bundle) becomes trivial. But this universal principal bundle and the associated finite-rank vector bundles remain nontrivial.
          \end{remark}

\end{enumerate}

\section{Summary and concluding remarks}

In this article, we have presented a unified treatment of Stiefel manifolds over Hilbert spaces, emphasizing the interplay between finite and infinite dimensions, real and complex scalars. In particular, we have shown that the polar decomposition provides a remarkably simple homeomorphism that allows a clean decomposition of homotopy types. \\

We conclude with a comprehensive table summarizing the connectivity properties of Stiefel manifolds in all dimensions and for both ground fields.

\begin{table}[h!]
    \centering
    \resizebox{\textwidth}{!}{
        \begin{tabular}{|l|c|c|c|c|}
            \hline
            \textbf{Manifold}                              & \textbf{Ground field} & \textbf{Dimension condition} & \textbf{Connected?} & \textbf{Higher homotopy}                                                                           \\
            \hline
            $\operatorname{St}_{\operatorname{orth}}(n,H)$ & $\mathbb{R}$          & $\dim H = n$                 & No (2 components)   & $\pi_0 \cong \mathbb{Z}_2$                                                                         \\
            \hline
            $\operatorname{St}_{\operatorname{orth}}(n,H)$ & $\mathbb{R}$          & $\dim H \geq n+1$, finite    & Yes                 & $(k-n-1)$-connected for $\dim H = k$                                                               \\
            \hline
            $\operatorname{St}_{\operatorname{orth}}(n,H)$ & $\mathbb{R}$          & $\dim H = \infty$            & Yes                 & Contractible                                                                                       \\
            \hline
            $\operatorname{St}_{\operatorname{orth}}(n,H)$ & $\mathbb{C}$          & $\dim H \geq n$, finite      & Yes                 & $(2k-2n)$-connected for $\dim H = k$                                                               \\
            \hline
            $\operatorname{St}_{\operatorname{orth}}(n,H)$ & $\mathbb{C}$          & $\dim H = \infty$            & Yes                 & Contractible                                                                                       \\
            \hline
            $\operatorname{St}(n,H)$                       & $\mathbb{R}$          & $\dim H = n$                 & No (2 components)   & $\pi_0 \cong \mathbb{Z}_2$                                                                         \\
            \hline
            $\operatorname{St}(n,H)$                       & $\mathbb{R}$          & $\dim H \geq n+1$, finite    & Yes                 & $\operatorname{St}(n,\mathbb{R}^k) \simeq \operatorname{St}_{\operatorname{orth}}(n,\mathbb{R}^k)$ \\
            \hline
            $\operatorname{St}(n,H)$                       & $\mathbb{R}$          & $\dim H = \infty$            & Yes                 & Contractible                                                                                       \\
            \hline
            $\operatorname{St}(n,H)$                       & $\mathbb{C}$          & $\dim H \geq n$, finite      & Yes                 & $\operatorname{St}(n,\mathbb{C}^k) \simeq \operatorname{St}_{\operatorname{orth}}(n,\mathbb{C}^k)$ \\
            \hline
            $\operatorname{St}(n,H)$                       & $\mathbb{C}$          & $\dim H = \infty$            & Yes                 & Contractible                                                                                       \\
            \hline
        \end{tabular}
    }
    \caption{Connectivity of Stiefel manifolds}
\end{table}

\begin{remark}
    The key points to remember:
    \begin{itemize}
        \item \textbf{Connectivity of $\operatorname{St}_{orth}(n,\mathbb{F}^k)$}:
              \begin{itemize}
                  \item For real scalars, the orthonormal Stiefel manifold has two components precisely when $k = n$; otherwise, for $k > n$, it is connected.
                  \item For complex scalars, the orthonormal Stiefel manifold is always connected when $k \geq n$.
              \end{itemize}
        \item \textbf{Homotopy type of $\operatorname{St}_{orth}(n,H)$ ($\dim H = \infty$)}: $\operatorname{St}_{orth}(n,\mathbb{F}^k)$ is contractible.
        \item The Stiefel manifold $\operatorname{St}(n,H)$ inherits the homotopy type from the orthonormal case because $\operatorname{P}_n(\mathbb{F})$ is contractible.
    \end{itemize}
\end{remark}

Several directions for further research remain open:

\begin{enumerate}
    \item \textbf{Equivariant theory:} The actions of compact groups on Stiefel manifolds and their equivariant homotopy types are only partially understood \cite{Bredon1972}.
    \item \textbf{Smooth structures:} While infinite-dimensional Stiefel manifolds are Hilbert manifolds \cite{Hamilton1982}, the existence of exotic smooth structures remains an open problem \cite{Milnor1964}.
    \item \textbf{Kuiper theorem for Fréchet spaces:} The validity of the Kuiper theorem for Fréchet spaces and its implications for infinite-dimensional Stiefel manifolds modeled on Fréchet spaces is open.
    \item \textbf{Stiefel manifolds over the quaternions:} Stiefel manifolds over quaternionic Hilbert spaces and their characteristic classes (Pontryagin classes for $\mathbb{H}$) have been less studied.
\end{enumerate}

We hope this exposition serves both as an introduction and as a reference for researchers working at the interface of algebraic topology and global analysis.

\section*{Acknowledgments}

This research has received no external funding.

\nocite{*} 
\bibliographystyle{amsplain}
\bibliography{references}

\Addresses

\end{document}